\documentclass[11pt]{article}
\usepackage[utf8]{inputenc}
\usepackage{comment}
\usepackage{amsmath}
\usepackage{amsthm}
\usepackage[bookmarks=false]{hyperref}
\usepackage{amssymb}
\usepackage{enumerate}
\usepackage{xcolor}
\usepackage{cases}
\usepackage{mathtools}
\usepackage{tikz} 
\DeclarePairedDelimiter\floor{\lfloor}{\rfloor}
\newtheorem{theo}{Theorem}[section]
\newtheorem{definition}{Definition}[section]
\newtheorem{prop}[theo]{Proposition}
\newtheorem{lemma}[theo]{Lemma}

\newtheorem{claim}[theo]{Claim}

\begin{document}
\date{}

\title{
An Asymptotically Sharp Bound on the Maximum Number of Independent Transversals
}

\author{Jake Ruotolo
\thanks
{Department of Mathematics, University of Central Florida, Orlando, FL 32816. \\Current Address: School of Engineering and Applied Sciences,
Harvard University, Cambridge,  MA 02134.
Email: {\tt jakeruotolo@g.harvard.edu}. 
}
\and 
Kevin Wang
\thanks{
Department of Mathematics \& Statistics, Grinnell College, Grinnell, IA 50112. \\Current Address: Department of Mathematics, University of Iowa, Iowa City, IA 52242.\\
Email: {\tt kevin-wang@uiowa.edu}.
}
\and
Fan Wei
\thanks
{
Department of Mathematics, Princeton University, Princeton, NJ
08544. \\
Current Address: Duke University, 120 Science Dr, Durham, NC 27710.
Email: {\tt fan.wei@duke.edu}.
Research supported by NSF Award
DMS-1953958.} 
}

\maketitle
\begin{abstract}
Let $G$ be a multipartite graph with partition $V_1, V_2,\ldots, V_k$ of $V(G)$. Let $d_{i,j}$ denote the edge density of the pair $(V_i, V_j)$. An independent transversal is an independent set of $G$ with exactly one vertex in each $V_i$. In this paper, we prove an asymptotically sharp upper bound on the maximum number of independent transversals given the $d_{i,j}$'s. 

\vspace{0.2cm}

\noindent
AMS Subject classification: 05C35, 05C69.  Keywords: independent transversal.
\end{abstract}

\section{Introduction}

Let $G$ be a multipartite graph with vertex partition $V_1, V_2, \ldots, V_k$. Jacob Fox asked the following question: given the edge density between every two vertex parts, what is the asymptotically maximum number of independent transversals in $G$ as $|V_i|$ goes to infinity for each $i$? An \emph{independent transversal} is an independent set of $G$ with exactly one vertex in each $V_i$. More precisely, for each $i \neq j$, let constant $d_{i,j}$ be the edge density between $V_i, V_j$, defined as $e(V_i, V_j)/|V_i||V_j|$. Independent transversals arise naturally in extremal combinatorics and bounding its number appears, for example, in inducibility type problems \cite{FHL2,FHL}. 
The previous known bound to this question is the following result of Fox, Huang, and Lee, which is an ingredient in \cite{FHL} to prove a bound on the number of induced copies of a given subgraph in another graph.

\begin{theo}[Lemma 4.1 \cite{FHL}] \label{thm:fox}
Let $k \geq 2$ be an integer. For each integer pair $1\leq i < j \leq k$, let $d_{i,j} = d_{j,i}$ be constants in $[0,1]$. Let $G$ be a multipartite graph with vertex partition $V_1, V_2, \ldots, V_k$ such that for each pair $1 \leq i < j \leq k$, the edge density between $V_i, V_j$ is $d_{i,j}$. 
Let $|V_i| = n_i$. 
Then the number of independent transversals in $G$ is at most
 \[\left(\prod_{1 \leq i < j \leq k}(1-d_{i,j})^{\lfloor k/2 \rfloor / \binom{k}{2}}\right)\left(\prod_{i=1}^k n_i\right).\]
\end{theo}
However, this bound is not sharp in general, or even not asymptotically sharp.  This means that there are choices of constants $d_{i,j}$'s such that the number of independent transversals divided by $\prod_{i=1}^k |V_i|$ is strictly less than $\prod_{1 \leq i < j \leq k}(1-d_{i,j})^{\lfloor k/2 \rfloor / \binom{k}{2}}$ as $ \prod_{i=1}^k |V_i|$ goes to infinity. 

 

In this paper, we prove an asymptotically sharp bound, previously asked by Jacob Fox \cite{pc}.  
Before stating our main theorem, we need the following definition. 
\begin{definition}\label{defn:odd cycle decomp}
An \emph{odd cycle decomposition} $H$ of the complete graph on $k$ vertices $K_k$ is a collection of disjoint multigraphs $F_1, F_2, \ldots, F_{\ell}$ satisfying $\bigcup_{i\in\ell}V(F_i)=V(K_k)$ such that for all $i\in\left[\ell\right]$, $F_i$ is an odd cycle, a double edge, or an isolated vertex. A double edge is obtained by adding an additional edge between the ends of an isolated edge. 
\end{definition}
\par
\begin{figure}\label{fig:odd cycle decomp}
\centering
\begin{tikzpicture}[node distance={20mm}, thick, main/.style = {draw, circle}]
\node[main] (1) {$v_1$}; 
\node[main] (2) [above right of=1] {$v_2$};
\node[main] (3) [right of=2] {$v_3$}; 
\node[main] (5) [below right of=1] {$v_5$}; 
\node[main] (4) [right of=5] {$v_4$};
\node[main] (6) [below of= 5] {$v_6$};
\node[main] (7) [below of= 4] {$v_7$};
\node[main] (8) [below right of= 3] {$v_8$};
\draw (1)--(2);
\draw (2)--(3);
\draw (3)--(4);
\draw (4)--(5);
\draw (1)--(5);
\path (6) edge [bend left] (7);
\path (6) edge [bend right] (7);
\end{tikzpicture} 
\caption{An odd cycle decomposition of $K_8$.}
\end{figure}
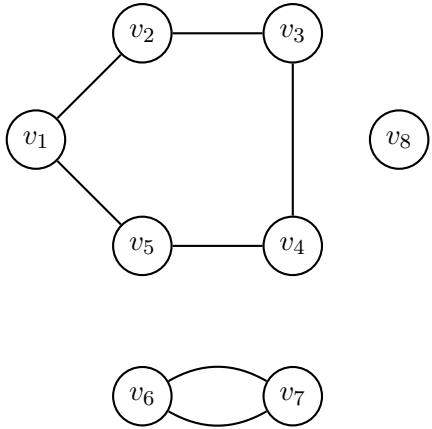

 Note that our definition of odd cycle decomposition may be different from other uses in the literature. See Figure~\ref{fig:odd cycle decomp} for an example. As a matter of notation, we denote an edge between vertices $v_i,v_j$ by $ij$.

\par Our main result is the following:

\begin{theo}\label{thm:main}
Let $k \geq 2$ be an integer. For each integer pair $1\leq i < j \leq k$, let $d_{i,j} = d_{j,i}$ be constants in $[0,1]$. Let $G$ be a multipartite graph with vertex partition $V_1, V_2, \ldots, V_k$ such that for each pair $1 \leq i < j \leq k$, the edge density between $V_i, V_j$ is $d_{i,j}$. 
Let $|V_i| = n_i$. Then the number of independent transversals in $G$ is at most \[ \min_{H \in \mathcal{H}} \left\{\prod_{F\in H} \prod_{ij\in E(F)}\sqrt{1-d_{i,j}}\right\}\prod_{i=1}^k n_i,\] 
where $\mathcal{H}$ is the set of all odd cycle decompositions of $K_k$. 
Furthermore, this bound is asymptotically sharp.
\end{theo}

By asymptotically sharp we mean that the bound is sharp up to a $o(n_1\cdots n_k)$-term, which refers to a function $f(n_1,\dots,n_k)$ with the property that\\ $\lim_{n_1,n_2,\ldots, n_k\to\infty}f(n_1,\dots,n_k)/(n_1\cdots n_k) = 0$. In particular, given densities $d_{i,j}$ for $i\neq j\in [k]$, we can construct a $k$-partite graph $G$ that attains the bound in Theorem~\ref{thm:main} up to a $o(n_1\cdots n_k)$-term. 

Observe that Theorem~\ref{thm:main} implies Theorem~\ref{thm:fox} as follows: Let $G$ be a multipartite graph as in the hypotheses of the two theorems. Let $\{a_tb_t\}_{t=1}^{\lfloor k/2 \rfloor}$ with $a_t,b_t\in[k]$ be a set of $\floor{k/2}$ edges corresponding to a matching on $K_k$. This corresponds to an odd cycle decomposition of $K_k$ where each subgraph $F_i$ is a double edge between $a_t$ and $b_t$. Since each edge $ij \in \binom{[k]}{2}$ is in $\lfloor k/2\rfloor / \binom{k}{2}$ fraction of all possible matchings as above, the bound in Theorem~\ref{thm:fox} is equal to the geometric mean over all products ${\prod_{t=1}^{\lfloor k/2 \rfloor} (1-d_{a_t, b_t})\prod_{i=1}^k n_i}$ given by all matchings. This geometric mean is at least the minimum over all possible products, which is at least the bound in Theorem~\ref{thm:main} on the maximum number of independent transversals in $G$. 

\section{Proof of Theorem~\ref{thm:main}}
The proof of Theorem~\ref{thm:main} comprises of two parts. In Section~\ref{subsec:ub}, we show that the result is an upper bound. In Section~\ref{subsec:as}, we show that the result is asymptotically sharp.
\subsection{Proof of Upper Bound}\label{subsec:ub}
The following lemma is a special case of Theorem~\ref{thm:main} (except that we consider transversal cliques instead of independent transversals). We use this case in the proof of Lemma~\ref{lem:cycle bound}, which is key to proving Theorem~\ref{thm:main}.

\begin{lemma}\label{lem:triangle}
For each integer pair $1\leq i<j\leq 3$, let $d_{i,j}=d_{j,i}$ be constants in $[0,1]$. 
Let $G$ be a tripartite graph with vertex partition $V_1, V_2, V_3$ such that for each pair $1\leq i<j\leq 3$, the edge density between $V_i,V_j$ is $d_{i,j}$. Let $\left|V_i\right| = n_i$. Suppose $d_{1,2}\leq d_{1,3} \leq d_{2,3}$. Then the number of transversal cliques in $G$ is at most \begin{equation}
    \min\left(d_{1,2}, \, \sqrt{d_{1,2}d_{1,3}d_{2,3}}\right)n_1n_2n_3. \label{eq:trianglebound}
\end{equation}
\end{lemma}
\begin{proof}
First suppose $d_{1,2} \leq \sqrt{d_{1,2}d_{1,3}d_{2,3}}$. Then we must prove that the number of transversal cliques in $G$ is at most $d_{1,2}n_1n_2n_3$. However, there are at most $d_{1,2}n_1n_2$ choices for an edge between $V_1$ and $V_2$, so clearly the number of transversal cliques is upper-bounded by $d_{1,2}n_1n_2n_3$. Thus we may assume $d_{1,2} > \sqrt{d_{1,2}d_{1,3}d_{2,3}}$. In particular, this implies $d_{1,2}, d_{1,3}, d_{2,3} > 0$. Let $V_1 = \{v_1, v_2, \ldots, v_{n_1}\}$. For each $i\in\left[n_1\right]$, we define $A_i = \{ v \in V_2 \mid v_iv \in E(G) \}$ and $B_i = \{ v \in V_3 \mid v_iv \in E(G) \}$. Let $\left|A_i\right| = a_in_2$, and $\left|B_i\right| = b_in_3$ for some $a_i, b_i \in [0,1]$. The number of transversal cliques in $G$ is $\sum_{i\in \left[n_1\right]}e(A_i, B_i)$. Since $0\leq e(A_i,B_i)\leq n_2n_3\min(d_{2,3}, a_ib_i)$, the number of transversal cliques in $G$ is at most $n_2n_3\sum_{i\in \left[n_1\right]}\min(d_{2,3}, a_ib_i)$. Note that we have the constraints $\sum_{i \in [n_1]} |A_i| = d_{1,2}n_1n_2$ and $\sum_{i \in [n_1]} |B_i| = d_{1,3}n_1n_3$. It suffices to solve the following problem:
\begin{align}
    &\text{Max } \sum_{i\in\left[n_1\right]}\min(a_ib_i, d_{2,3}) \label{eq:obj} \\
    &\text{subject to }
    \sum_{i\in \left[n_1\right]}a_i \leq  d_{1,2}n_1, \sum_{i\in \left[n_1\right]}b_i \leq d_{1,3}n_1, \text{ where } 0 \leq a_i, b_i \leq 1. \label{eq:constraint}
\end{align}
Here we still have the assumption that $d_{1,2} \leq d_{1,3}$. 

\par Let ${\mathcal{O}}$ be the set of all optimal solutions to our problem. The set $\mathcal{O}$ is nonempty because the objective function (\ref{eq:obj}) is continuous on the compact set defined by
the constraints (\ref{eq:constraint}). We prove a series of claims, using local adjustments, to show that there must be optimal solutions in $\mathcal{O}$ that satisfy certain nice conditions. We then use these conditions to get the desired upper bound on the number of transversal cliques in $G$.

\begin{claim}\label{claim:01}
There is a nonempty subset $\mathcal{O}_1$ of $\mathcal{O}$ such that for each \\$O=\{a_1,\dots,a_{n_1},b_1,\dots,b_{n_1}\} \in \mathcal{O}_1$, we have $a_ib_i \leq d_{2,3}$ for all $i$.
\end{claim}
\begin{proof}
Suppose that there is $O \in \mathcal{O}$ where for some $i \in [n_1]$, $a_ib_i > d_{2,3}$. Then we can decrease $a_i, b_i$ to $a_i'$, $b_i'$ such that $a_i' b_i' = d_{2,3}$. This new set of variables still satisfies the constraints (\ref{eq:constraint}) while the objective function (\ref{eq:obj}) has the same value. 
\end{proof}

\begin{claim}\label{claim:02}
Let $O=\{a_1,\dots,a_{n_1},b_1,\dots,b_{n_1}\}\in\mathcal{O}_1$. Then $\sum_{i\in[n_1]}a_i\leq \sum_{i\in[n_1]}b_i$.
\end{claim}
\begin{proof}
Assume, to the contrary, that $\sum_{i\in[n_1]}a_i > \sum_{i\in[n_1]}b_i$. This implies $\sum_{i\in[n_1]}b_i < d_{1,3}n_1$. To obtain a contradiction, we have the following two cases (note that Claim~\ref{claim:01} implies $a_ib_i\leq d_{2,3}$ for all $i$).
\par First assume that there exists an $i\in[n_1]$ such that $a_ib_i < d_{2,3}$ and $a_i > b_i$. Then we can choose $\varepsilon > 0$ such that $\varepsilon \leq \min(d_{1,3}n_1 - \sum_{i\in[n_1]}b_i, 1-b_i)$. By increasing $b_i$ by $\varepsilon$, we satisfy the constraints (\ref{eq:constraint}) and strictly increase the objective function (\ref{eq:obj}), contradicting the optimality of $O$. This implies that for all $i\in[n_1]$ with $a_ib_i < d_{2,3}$ we have $a_i\leq b_i$. 
\par Now assume that there exists $i\in[n_1]$ such that $a_ib_i = d_{2,3}$ and $a_i > b_i$. We prove that there exists $j\in[n_1]$ with $a_jb_j<d_{2,3}$. To this end, assume for contradiction that for all $j\in[n_1]$, we have $a_jb_j=d_{2,3}$. Then $d_{2,3}n_1=\sum_{j\in[n_1]}a_jb_j$. Since $a_i>b_i$, we have $a_i(1-b_i)>0$, and so $\sum_{j\in[n_1]}a_j-a_jb_j=\sum_{j\in[n_1]}a_j(1-b_j)>0$. Thus $\sum_{j\in[n_1]}a_j>\sum_{j\in[n_1]}a_jb_j$, which implies $d_{1,2}>d_{2,3}$, a contradiction. Therefore we can fix $j\in[n_1]$ with $a_jb_j < d_{2,3}$. 
By the preceding paragraph, we must have $a_j<1$. Choose $\varepsilon>0$ such that $\varepsilon\leq\min(a_i-b_i,1-a_j, d_{1,3}n_1 -\sum_{i\in[n_1]}b_i)$. Then $(a_i-\varepsilon)(b_i+\varepsilon) \geq a_ib_i$ and $ (a_j + \varepsilon)b_j > a_jb_j$. Hence, we can increase the objective function while satisfying the constraints, once again contradicting the optimality of $O$. Thus for all $i\in[n_1]$ such that $a_ib_i=d_{2,3}$, we have $a_i\leq b_i$.
\par The two cases imply $\sum_{i\in[n_1]}a_i\leq \sum_{i\in[n_1]}b_i$, contradicting our initial assumption.
\end{proof}

\begin{claim}\label{claim:03}
Let $O=\{a_1,\dots,a_{n_1},b_1,\dots,b_{n_1}\}\in\mathcal{O}_1$. For all $i \in \left[n_1\right]$ with $a_ib_i < d_{2,3}$, we have $a_i \leq b_i$.
\end{claim}
\begin{proof}
Assume for contradiction that there exists $i \in \left[n_1\right]$ with $a_ib_i < d_{2,3}$ but $a_i > b_i$. Then by Claim~\ref{claim:02}, there must exist $j \in \left[n_1\right]$ with $a_j < b_j$. Notice that if $\varepsilon \in (0, \min(a_i - b_i, b_j - a_j))$, then $(a_i - \varepsilon)(b_i + \varepsilon) - a_ib_i = \varepsilon(a_i - b_i) - \varepsilon^2 > 0$ and $(a_j + \varepsilon)(b_j - \varepsilon) - a_jb_j = \varepsilon(b_j-a_j)-\varepsilon^2 > 0$. Since $a_ib_i < d_{2,3}$, we can choose such an $\varepsilon$ to obtain a solution to the objective function whose value is greater than that of $O$ as long as the constraints are still satisfied. Thus, by choosing $\varepsilon \in (0, \min(a_i - b_i, b_j - a_j, 1-b_i, a_i, b_j, 1-a_j))$, we increase the objective function while satisfying the constraints. This contradicts the optimality of~$O$.
\end{proof}

Let $\mathcal{O}_2$ be the set of $O$ in $\mathcal{O}_1$ where $\sum a_i + \sum b_i$ is minimized. 
\begin{claim}\label{claim:equality}
Let $O = \{a_1, \ldots, a_{n_1}, b_1,\ldots, b_{n_1}\} \in\mathcal{O}_2$. For all $i,j\in\left[n_1\right]$ with $a_ib_i = a_jb_j = d_{2,3}$, we have $a_i = a_j$ and $b_i = b_j$.
\end{claim}
\begin{proof}
Assume for contradiction that there exists $i, j\in[n_1]$ with $a_ib_i = a_jb_j = d_{2,3}$ but $a_i \neq a_j$. By symmetry, say $a_i > a_j$. Then $b_i < b_j$. Choose sufficiently small $\varepsilon, \varepsilon' > 0$ and let $a_i' = a_i - \varepsilon, a_j' = a_j+\varepsilon, b_i' = b_i + \varepsilon', b_j' = b_j - \varepsilon'$. Now $a_i' b_i'  + a_j' b_j'= (a_i - \varepsilon)(b_i + \varepsilon') + (a_j + \varepsilon)(b_j - \varepsilon') = 2d_{2,3} + \varepsilon'(a_i - a_j) + \varepsilon(b_j - b_i) - 2\varepsilon\varepsilon' > 2d_{2,3}$ when $\varepsilon, \varepsilon'$ are sufficiently small. Observe $a_i'b_i' - a_j'b_j' = \varepsilon'(a_i + a_j) - \varepsilon(b_i + b_j)$. Thus if $\varepsilon'/\varepsilon = (b_i + b_j)/(a_i+ a_j)$ with $\varepsilon,\varepsilon'$ sufficiently small, then we must have $a_i' b_i' =a_j' b_j' > d_{2,3}$. Therefore we can decrease the values of $a_i', a_j'$ so that $a_i'b_i' = a_j'b_j' = d_{2,3}$, and in this way the objective function value does not change although $\sum a_i + \sum b_i$
strictly decreases. This contradicts the definition of $\mathcal{O}_2$.
\end{proof}

For each $O = \{a_1, \ldots, a_{n_1}, b_1,\ldots, b_{n_1}\} \in\mathcal{O}_2$, define $C_O = \{ i \in [n_1] \mid 0 < a_ib_i < d_{2,3} \}$ and $D_O = \{ i \in [n_1] \mid a_ib_i=d_{2,3} \}$. We prove the following two claims involving $C_O$. 

\begin{claim}\label{claim:b_i < 1}
Let $O = \{a_1, \ldots, a_{n_1}, b_1,\ldots, b_{n_1}\} \in\mathcal{O}_2$. Either the objective function is at most $\sqrt{d_{1,2}d_{1,3}d_{2,3}} n_1$ (and so the desired upper bound (\ref{eq:trianglebound}) holds), or for all $i \in C_O$, we have $b_i < 1$.
\end{claim}
\begin{proof}
Suppose that there exists $i\in C_O$ with $b_i=1$. Then we must prove that the objective function is at most $\sqrt{d_{1,2}d_{1,3}d_{2,3}}n_1$. Since $b_i=1$, the definition of $C_O$ implies $a_i<1$. 
Assume for contradiction that there exists $j$ with $a_jb_j>0$ and $b_j<1$. Then we can decrease $a_j$ and increase $a_i$ by some sufficiently small $\varepsilon > 0$ to contradict the optimality of $O$.
Thus for all $j$ with $a_j b_j > 0$, we must have $b_j=1$. Since $\sum_{i\in[n_1]}b_i\leq d_{1,3}n_1$ (recall constraints (\ref{eq:constraint})), the number of $j$ such that $a_jb_j > 0$ is at most $d_{1,3} n_1$. 
The objective function is $\sum_{j \in [n_1]} \min(a_j b_j, d_{2,3})= \sum_{j \in [n_1], a_jb_j > 0} \min(a_jb_j, d_{2,3}) \leq d_{1,3} n_1 d_{2,3}$. By our assumption that $d_{1,2} > \sqrt{d_{1,2}d_{1,3}d_{2,3}}$, we have $d_{1,3}d_{2,3}<\sqrt{d_{1,2}d_{1,3}d_{2,3}}$, and so $\sum_{j \in [n_1]} \min(a_j b_j, d_{2,3}) < \sqrt{d_{1,2}d_{1,3}d_{2,3}} n_1$, as desired. 
\end{proof}

\begin{claim}\label{claim:C1}
Let $O = \{a_1, \ldots, a_{n_1}, b_1,\ldots, b_{n_1}\} \in\mathcal{O}_2$. Then 
$|C_O| \leq 1$. 
\end{claim}
\begin{proof}
\par Assume for contradiction that $|C_O| \geq 2$. Without loss of generality, assume indices $1,2 \in C_O$. By definition of $C_O$ and Claim~\ref{claim:03}, we have $0 < a_1, a_2 < 1$. We have three cases. If $b_1<b_2$, then we can decrease $a_1$ and increase $a_2$ by some $\varepsilon\in(0,\min(a_1,1-a_2)]$ to increase the objective function, thereby contradicting the optimality of $O$. Similarly, if $b_1>b_2$, then we can decrease $a_2$ and increase $a_1$ by some $\varepsilon\in(0,\min(a_2,1-a_1)]$ to contradict the optimality of $O$. If $b_1=b_2$, then we can assume $0<b_1 = b_2 < 1$ by Claim~\ref{claim:b_i < 1} and the definition of $C_O$. Now, without loss of generality, suppose $a_1 \leq a_2$. In this case, we can decrease $a_1,b_1$ and increase $a_2,b_2$ by the same amount $\varepsilon\in(0,\min(a_1,1-a_2, b_1, 1-b_2)]$ to contradict the optimality of $O$. 
\end{proof}

\par To finish the proof, fix $O=\{a_1,\dots,a_{n_1},b_1,\dots,b_{n_1}\}\in\mathcal{O}_2$. 
Claim~\ref{claim:equality} implies that for each $i, i' \in D_O$, we have $a_i = a_{i'} = a$ and $b_i = b_{i'} = b$ for some $a,b \in (0,1]$ with $ab=d_{2,3}$. The objective function is \[ \sum_{i\in [n_1]}\min(a_ib_i, d_{2,3}) = \sum_{i\in C_O\cup D_O}\min(a_ib_i, d_{2,3}) = d_{2,3}|D_O| + \sum_{j\in C_O}a_jb_j.\] By the constraints, the quantity $|D_O|$ satisfies $|D_O| a + \sum_{j \in C_O} a_j \leq d_{1,2}n_1$ and $|D_O| b + \sum_{j \in C_O}b_j \leq d_{1,3}n_1$. Using this, we can write
\begin{align*}
 |D_O| &\leq \sqrt{(d_{1,2}n_1-\sum_{j\in C_O}a_j)(d_{1,3}n_1-\sum_{j\in C_O}b_j)/d_{2,3}}.
\end{align*}
If $\left| C_O\right| = 0$, then $\left|D_O\right| \leq n_1\sqrt{d_{1,2}d_{1,3}/d_{2,3}}$. This implies that the objective function is at most $n_1\sqrt{d_{1,2}d_{1,3}d_{2,3}}$, and so the number of transversal cliques in $G$ is at most $n_1n_2n_3\sqrt{d_{1,2}d_{1,3}d_{2,3}}$.

By Claim~\ref{claim:C1}, the only other case that we need to consider is $\left|C_O\right| = 1$. Let $j\in C_O$. In this case, we can bound $|D_O|$ by
\begin{align*}
    \left|D_O\right|&\leq \sqrt{(d_{1,2}n_1-a_j)(d_{1,3}n_1-b_j)/d_{2,3}}\\ 
    &=\sqrt{1/d_{2,3}}\sqrt{d_{1,2}d_{1,3}n^2_1-(d_{1,2}b_j+d_{1,3}a_j)n_1+a_jb_j}\\
    &\leq \sqrt{1/d_{2,3}}\sqrt{d_{1,2}d_{1,3}n^2_1-2\sqrt{d_{1,2}d_{1,3}a_jb_j}n_1+a_jb_j}\\
    &= \sqrt{1/d_{2,3}}\left(\sqrt{d_{1,2}d_{1,3}n^2_1} - \sqrt{a_jb_j}\right),
\end{align*}
where the second inequality follows from the AM-GM inequality ($x+y \geq 2\sqrt{xy}$ for any $x,y\geq 0$). In the last step above, note that $\sqrt{d_{1,2}d_{1,3}n^2_1} - \sqrt{a_jb_j}\geq 0$ because $a_j\leq d_{1,2}n_1$ and $b_j\leq d_{1,3}n_1$ by the constraints. Using the fact that $0<a_jb_j<d_{2,3}$, we can bound the objective function by
\begin{align*}
    d_{2,3}\left|D_O\right| + a_jb_j &\leq  n_1\sqrt{d_{1,2}d_{1,3}d_{2,3}} - \sqrt{a_jb_jd_{2,3}} + a_jb_j\\
    &< n_1\sqrt{d_{1,2}d_{1,3}d_{2,3}}.
\end{align*}
This implies that the number of transversal cliques in $G$ is at most $n_1n_2n_3\sqrt{d_{1,2}d_{1,3}d_{2,3}}$, completing the proof.

\end{proof}

We use Lemma~\ref{lem:triangle} in the below proof of Lemma~\ref{lem:cycle bound}, which is key to the proof of the upper bound in Theorem~\ref{thm:main}. 

\begin{lemma}\label{lem:cycle bound}
Let $k \geq 3$ be an integer. For each integer pair $1\leq i<j\leq k$, let $d_{i,j}=d_{j,i}$ be constants in $[0,1]$. Let $G$ be a multipartite graph with vertex partition $V_1, V_2, \ldots, V_{k}$ such that for each pair $1\leq i<j\leq k$, the edge density between $V_i,V_j$ is $d_{i,j}$. Let $|V_i| = n_i$. Then the number of transversal cliques in $G$ is at most
\begin{align*}
\left(\prod_{i = 1}^{k} \sqrt{d_{i,i+1}}\right)\left(\prod_{i=1}^k n_i\right),
\end{align*}
where the index $i+1$ is modulo $k$.
\end{lemma}

\begin{proof}
We show by induction on $k$ that the statement holds for all $k\geq 3$.
By Lemma~\ref{lem:triangle}, the statement holds for $k=3$. Assume that the statement holds for $k-1 \geq 3$.
Let $N(v) = \{ u : uv \in E(G) \}$. For $v_i \in V_1$, let $A^i_{1,2} = N(v_i)\cap V_2$ and $A^i_{1,k} = N(v_i)\cap V_k$. Let $\left|A^i_{1,2}\right| = a^i_{1,2}n_2$ and $\left|A^i_{1,k}\right| = a^i_{1,k}n_k$. We bound the number of transversal cliques in $G$ containing $v_i$ (i.e., the number of transversal cliques in $G[A^i_{1,2},V_3, \ldots, V_{k-1}, A^i_{1,k}]$), and then sum this quantity over all $v_i\in V_1$ to get the bound in the lemma. 
\par If $a^i_{1,2} = 0$ or $a^i_{1,k} = 0$, then there are zero transversal cliques in $G$ containing $v_i$, so we need only consider $i$ with $a^i_{1,2}, a^i_{1,k} > 0$. Observe that the edge density between $A^i_{1,2},V_3$ is at most $\min(d_{2,3}/a^i_{1,2},1)$, the edge density between $A^i_{1,k},V_{k-1}$ is at most $\min(d_{k-1,k}/a^i_{1,k},1)$, and the edge density between $A^i_{1,2},A^i_{1,k}$ is at most 1. Then applying the inductive hypothesis to $G[A^i_{1,2},V_3, \ldots, V_{k-1}, A^i_{1,k}]$, the number of transversal cliques in $G[A^i_{1,2},V_3, \ldots, V_{k-1}, A^i_{1,k}]$ is at most 
\begin{align*}
   a^i_{1,2}n_2n_3\cdots n_{k-1}a^i_{1,k}n_k\sqrt{\min(d_{2,3}/a^i_{1,2}, 1)d_{3,4}\cdots d_{k-2, k-1}\min(d_{k-1,k}/a^i_{1,k}, 1)}.
\end{align*}
To bound the number of transversal cliques in $G$, we sum the above quantity over all $v_i\in V_1$. This amounts to bounding 
 \begin{align*}
    \sum_{v_i\in V_1}a^i_{1,2}a^i_{1,k}\sqrt{\min(d_{2,3}/a^i_{1,2},1)\min(d_{k-1,k}/a^i_{1,k}, 1)}
 \end{align*} subject to $\sum_{v_i\in V_1}a^i_{1,2} = d_{1,2}n_1$ and $\sum_{v_i\in V_1}a^i_{1,k} = d_{1,k}n_1$. In turn, this amounts to bounding
 \begin{align*}
      \sqrt{d_{2,3}d_{k-1,k}}\sum_{v_i\in V_1}\sqrt{a^i_{1,2}a^i_{1,k}}
 \end{align*}
 subject to $\sum_{v_i\in V_1}a^i_{1,2} = d_{1,2}n_1$ and $\sum_{v_i\in V_1}a^i_{1,k} = d_{1,k}n_1$. Applying the Cauchy-Schwarz inequality, we have
  \begin{align*}
      \sum_{v_i\in V_1}\sqrt{a^i_{1,2}a^i_{1,k}}&\leq  \sqrt{\sum_{v_i\in V_1}a^i_{1,2}\sum_{v_i\in V_1}a^i_{1,k}}\\
      &= n_1\sqrt{d_{1,2}d_{1,k}}.
 \end{align*}
 Hence, the number of transversal cliques in $G$ is at most
 \begin{align*}
      n_1n_2\cdots n_k\sqrt{d_{1,2}d_{2,3}\cdots d_{k-1,k}d_{k,1}},
 \end{align*}
 as desired.
\end{proof}

We now use Lemma~\ref{lem:cycle bound} to prove the below proposition, which is the upper bound in Theorem~\ref{thm:main}. Although we may minimize over all (not necessarily odd) cycle decompositions in the below bound, it turns out that it suffices to minimize over just the odd ones to get that the bound is asymptotically sharp (see Section~\ref{subsec:as}).

\begin{prop}\label{upper bound}
Let $k \geq 2$ be an integer. For each integer pair $1\leq i < j \leq k$, let $d_{i,j} = d_{j,i}$ be constants in $[0,1]$. Let $G$ be a multipartite graph with vertex partition $V_1, V_2, \ldots, V_k$ such that for each pair $1 \leq i < j \leq k$, the edge density between $V_i, V_j$ is $d_{i,j}$. 
Let $|V_i| = n_i$. Then the number of independent transversals in $G$ is at most \[ \min_{H \in \mathcal{H}} \left\{ \prod_{F\in H} \prod_{ij\in E(F)}\sqrt{1-d_{i,j}}\right\}\prod_{i=1}^kn_i,\] 
where $\mathcal{H}$ is the set of all odd cycle decompositions of $K_k$. 
\end{prop}

\begin{proof}
Since the number of independent transversals in $G$ is equal to the number of transversal cliques in $\overline{G}$, we obtain an upper bound for the latter. Note that each edge density $d_{i,j}$ in $G$ corresponds to the edge density $1-d_{i,j}$ in $\overline{G}$.
Let $H = \{F_1, F_2, \ldots, F_l\} \in \mathcal{H}$. By definition of odd cycle decomposition (Definition \ref{defn:odd cycle decomp}), the sets $V(F_1), V(F_2), \ldots, V(F_l)$ partition $V(K_k)$. For each $i \in [l]$, let $S_i = \{V_j: v_j \in V(F_i)\}$. Observe that $S_1, \dots, S_l$ correspond bijectively to $V(F_1), \dots, V(F_l)$ and partitions $V_1, V_2, \ldots, V_k$. 
We now obtain a bound for each possible $\overline{G}\left[S_i\right]$, that is, each possible subgraph of $\overline{G}$ induced by $S_i$. If $S_r = \{V_1, \dots, V_m\}$ corresponds to an odd cycle, then the number of transversal cliques in $\overline{G}\left[S_r\right]$ is at most $n_1 \cdots n_m\prod_{ij\in E(F_r)}\sqrt{1-d_{i,j}}$ by Lemma~\ref{lem:cycle bound}. If $S_r = \{V_1, V_2\}$ corresponds to a double edge, then the number of transversal cliques in $\overline{G}\left[S_r\right]$ is $n_1n_2(1-d_{1,2}) = n_1n_2\prod_{ij\in E(F_r)}\sqrt{1-d_{i,j}}$. If $S_r = \{ V_1 \}$ corresponds to an isolated vertex, then the number of transversal cliques in $\overline{G}[S_r]$ is trivially $n_1$. Putting the three cases together, the number of transversal cliques in $\overline{G}$ is at most $\displaystyle n_1\cdots n_k \prod_{F\in H} \prod_{ij\in E(F)}\sqrt{1-d_{i,j}}$. Minimizing this over $\mathcal{H}$, we obtain the desired upper bound.
\end{proof}

The following section is dedicated to showing that the above upper bound is asymptotically sharp.

\subsection{Proof of Asymptotic Sharpness}\label{subsec:as}
Assume the hypotheses of Theorem~\ref{thm:main}.
To show that the upper bound in Theorem~\ref{thm:main} is asymptotically sharp, we use linear programming and duality to prove the existence of a $k$-partite graph $G$ that attains the bound up to a $o(n_1\cdots n_k)$-term. Since the number of independent transversals in $G$ equals the number of transversal cliques in $\overline{G}$, it suffices to consider the latter. Note that each edge density $d_{i,j}$ in $G$ corresponds to the edge density $1-d_{i,j}$ in $\overline{G}$. Consider the following construction of $\overline{G}$ which satisfies our constraints:
\begin{align*}
    &\text{For each }1 \leq i \leq k,\text{ let }S_i \subseteq V_i \text{ with }|S_i| = \floor{a_in_i}, \text{ where }\\&a_ia_j \leq 1-d_{i,j}\:(\text{for all } 1 \leq i < j \leq k), \text{ and }a_i\in(0,1]\:(\text{for all }i).\\&\text{Suppose that each pair }S_i, S_j\,(i \neq j)\,\text{ forms a complete bipartite graph}.
\end{align*}
In this construction, the number of transversal cliques is at least $\prod_{i=1}^k \floor{a_in_i}$. We will show that some construction described above is close to the maximum number of transversal cliques.
This motivates the following optimization problem:
\begin{align*}\label{lp: 1}
    \texttt{P}\text{: }&\text{Max } \prod_{i=1}^k a_in_i\\
    &\text{subject to } a_ia_j \leq 1-d_{i,j}\:(\text{for all } 1 \leq i < j \leq k), \, a_i\in(0,1]\:(\text{for all }i).
\end{align*}
If we have an optimal solution $(a_1,\ldots,a_k)$, then $\prod_{i=1}^k \floor{a_in_i}$ is a lower bound on the number of transversal cliques in the above construction in $\overline{G}$. 
Observe $\prod_{i=1}^k \floor{a_in_i}=\left(\prod_{i=1}^k a_in_i\right) + o(n_1\cdots n_k)$. Thus the construction will attain the bound in Theorem~\ref{thm:main} up to a $o(n_1\cdots n_k)$-term if we can prove that $\prod_{i=1}^k a_in_i$ is given by a product over some odd cycle decomposition of $K_k$ (see Equation~(\ref{eq:oddcycledecomp}) for more specificity). 

\par In problem \texttt{P}, note that $a_i>0$ for all $i$ implies $d_{i,j}<1$ for all $i,j$ (if $d_{i,j}=1$ is allowed, then
the number of transversal cliques in $\overline{G}$ is zero).
Applying the natural log to \texttt{P}, we get an equivalent linear programming problem:
\begin{align*}
    \texttt{LP1}\text{: }&\text{Max } \sum_{i=1}^k \ln(n_i) + \sum_{i=1}^k \ln(a_i)\\
    &\text{subject to }\ln(a_i) + \ln(a_j) \leq \ln(1-d_{i,j}),\, \ln(a_i) \leq 0.
\end{align*}
Ignoring the constant $\sum_{i=1}^k \ln(n_i)$ and setting $b_i = -\ln(a_i)$, $p_{i,j} = -\ln(1-d_{i,j})$, we can rewrite \texttt{LP1} as an equivalent problem more convenient to work with:
\begin{align*}
    \texttt{LP2}\text{: }&\text{Min } \sum_{i=1}^k b_i  \\
    &\text{subject to }b_i + b_j \geq p_{i,j}, \, b_i \geq 0.
\end{align*}
The dual of \texttt{LP2} is
\begin{align*}
    \texttt{LP2-dual}\text{: }&\text{Max } \sum_{1 \leq i < j \leq k} p_{i,j}x_{i,j}\\
    &\text{subject to }\sum_{i \neq j} x_{i,j} \leq 1 \text{ for each fixed }i\text{ }(1 \leq i \leq k), \, x_{i,j} \geq 0.
\end{align*}
Introducing the slack variable $y_i$, \texttt{LP2-dual} becomes 
\begin{align*}
    \texttt{LP2-dual$'$}\text{: }&\text{Max } \sum_{1 \leq i < j \leq k} p_{i,j}x_{i,j}\\
    &\text{subject to }y_i + \sum_{i \neq j} x_{i,j} = 1 \text{ for each fixed }i\text{ }(1 \leq i \leq k), x_{i,j} \geq 0, y_i \geq 0.
\end{align*}

We need the following results from linear programming. The first is known as the duality of linear programming.

\begin{theo}[Duality of Linear Programming~\cite{LP}]\label{duality}
For the linear programs
\begin{center}
maximize $\textbf{c}^{T}\textbf{x}$ subject to $\textbf{Ax}\leq \textbf{b}$ and $\textbf{x}\geq \textbf{0}$\hspace{1cm}(P)\\
\end{center}
and
\begin{center}
minimize $\textbf{b}^{T}\textbf{y}$ subject to $\textbf{A}^{T}\textbf{y}\geq \textbf{c}$ and $\textbf{y}\geq \textbf{0}$\hspace{1cm}(D),
\end{center} the following holds:
   If both (P) and (D) have a feasible solution, then both have an optimal solution, and if $\textbf{x}^*$ is an optimal solution of (P) and $\textbf{y}^*$ is an optimal solution of (D), then
    \begin{center}
        $\textbf{c}^T\textbf{x}^* = \textbf{b}^T\textbf{y}^*$.
    \end{center}
That is, the maximum of (P) equals the minimum of (D).
\end{theo}

We also need a corollary to the above duality theorem called complementary slackness.

\begin{theo}[Complementary Slackness~\cite{LP}]\label{compSlack}
Let $\textbf{x}^* = (x_1^*, x_2^*, \ldots, x_n^*)$ be a feasible solution of the linear program
\begin{center}
maximize $\textbf{c}^{T}\textbf{x}$ subject to $\textbf{Ax}\leq \textbf{b}$ and $\textbf{x}\geq \textbf{0}$\hspace{1cm}(P),\\
\end{center}
and let $\textbf{y}^* = (y_1^*, y_2^*, \ldots, y_m^*)$ be a feasible solution of the dual linear program
\begin{center}
  minimize $\textbf{b}^{T}\textbf{y}$ subject to $\textbf{A}^{T}\textbf{y}\geq \textbf{c}$ and $\textbf{y}\geq \textbf{0}$\hspace{1cm}(D).  
\end{center}
Then the following two statements are equivalent:
\begin{enumerate}
    \item $\textbf{x}^*$ is optimal for (P) and $\textbf{y}^*$ is optimal for (D).
    \item For all $i\in\left[m\right]$, $\textbf{x}^*$ satisfies the $i$th constraint of (P) with equality or $y_i^* = 0$; similarly, for all $j\in\left[n\right]$, $\textbf{y}^*$ satisfies the $j$th constraint of (D) with equality or $x_j^* = 0$.
\end{enumerate}
\end{theo}

Since \texttt{LP2-dual$'$} and \texttt{LP2} are both feasible, both have an optimal solution by Theorem~\ref{duality}. We will show that there is an optimal solution to \texttt{LP2-dual$'$} which is given by an odd cycle decomposition of $K_k$. Then Theorem~\ref{compSlack} will show that optimal solutions to \texttt{LP2} are given by that odd cycle decomposition of $K_k$, which will imply that $\prod_{i=1}^k a_in_i$ is given by a product over that odd cycle decomposition of $K_k$, as desired.

More specifically, we define the graph $Q$ on $k$ vertices $x_1, \ldots, x_k$ as follows (recall that $\texttt{LP2-dual}'$ has variables $x_{i,j}$): 
\begin{enumerate}
    \item there exists an edge $x_ix_j \in Q$ if and only if $x_{i,j} > 0$,
    \item there exists a self-loop on vertex $x_i$ if and only if $y_i > 0$.
\end{enumerate}
Thus $Q$ can be used to represent elements of the feasible set of \texttt{LP2-dual$'$} (non-edges and non-self-loops correspond to variables being $0$). Consider the set of $Q$'s that represent optimal solutions to \texttt{LP2-dual$'$} (this set is nonempty by Theorem~\ref{duality}). From this set, consider the $Q$'s with the minimum number of non-loop edges, and then among which choose the ones with the minimum number of self-loops. Call the resulting set $\mathcal{Q}$. Asymptotic sharpness of the bound in Theorem~\ref{thm:main} is proved in the following proposition. 

\begin{prop}\label{asymptotic sharpness}
Let $k \geq 2$ be an integer. For each integer pair $1\leq i < j \leq k$, let $d_{i,j} = d_{j,i}$ be constants in $[0,1]$. For any sufficiently large integers $n_1,\ldots, n_k$, there exists a multipartite graph $G$ with vertex partition $V_1, \ldots, V_k$ where $|V_i| = n_i$ for each $i\in [k]$ such that the number of independent transversals in $G$ is at least \[ \min_{H \in \mathcal{H}} \left\{\prod_{F\in H} \prod_{ij\in E(F)}\sqrt{1-d_{i,j}}\right\}\left(\prod_{i=1}^kn_i\right)+o\left(\prod_{i=1}^kn_i\right),\] 
where $\mathcal{H}$ is the set of all odd cycle decompositions of $K_k$, and the edge density between $V_i, V_j$ is $d_{i,j}$ for each pair $1 \leq i < j \leq k$. Hence, the bound in Theorem~\ref{thm:main} is asymptotically sharp. 
\end{prop}

\begin{proof}
Recall that we must prove that $\prod_{i=1}^k a_in_i$ is given by a product over some odd cycle decomposition of $K_k$. We show that Claim~\ref{claim:Q'} below gives the required odd cycle decomposition. We then prove Claim~\ref{claim:Q'} to complete the proof of the proposition.

\begin{claim}\label{claim:Q'}
For each $Q\in\mathcal{Q}$, define a graph $Q'$ constructed by adding an additional edge between the ends of isolated edges in $Q$ and removing self-loops from isolated vertices in $Q$. Then there is a $Q'$ constructed from some $Q \in \mathcal{Q}$ which is an odd cycle decomposition of $K_k$.
\end{claim}

Suppose that Claim~\ref{claim:Q'} holds. Fix such a $Q$ and $Q'$. Recall that $Q$ represents an optimal solution to \texttt{LP2-dual$'$}. We are done if we can prove 
\begin{equation}\label{eq:oddcycledecomp}
    \prod_{i=1}^k a_in_i = n_1\cdots n_k\prod_{F\in Q'} \prod_{ij\in E(F)}\sqrt{1-d_{i,j}}.
\end{equation}
Let $(b_1, \dots, b_k)$ be an optimal solution to \texttt{LP2}. 
Then Theorem~\ref{compSlack} (Complementary Slackness) and the construction of $Q$ imply that each edge $x_ix_j$ in $Q$ corresponds to an equality in \texttt{LP2}'s constraints, that is, $b_i + b_j = p_{i,j}$. 
Similarly, each isolated vertex $x_i$ of $Q$ corresponds to $b_i = 0$. Converting $(b_1, \dots, b_k)$ back to an optimal solution $(a_1, \dots, a_k)$ of \texttt{P}, we get that for each edge $x_ix_j$ in $Q$, we have an equality in \texttt{P}'s constraints, that is, $a_ia_j = 1-d_{i,j}$; 
similarly, for each isolated vertex $x_i$ of $Q$, we have the equation $a_i = 1$.
Thus, since $Q'$ is an odd cycle decomposition of $K_k$, the objective function of \texttt{P}, $\prod a_i n_i$, has the above form, as desired.

To prove Claim~\ref{claim:Q'}, we first need a series of claims.

\begin{claim}\label{self-loop}
For any $Q\in \mathcal{Q}$, there is at most one vertex with a self-loop in $Q$.
\end{claim}
\begin{proof}
Assume for contradiction that there exists $Q\in\mathcal{Q}$ with at least two vertices with a self-loop. Then by the construction of $Q$, there exist $i < j$ with $y_i > 0$ and $y_j > 0$. To contradict the optimality of $Q$, set $\varepsilon = \min(y_i, y_j)$. Let $y_i' = y_i - \varepsilon$, $y_j' = y_j - \varepsilon$, $x_{i,j}' = x_{i,j} + \varepsilon$, and all other primed variables have the same value as the corresponding non-primed variable. Then the objective function of \texttt{LP2-dual$'$}, $\sum p_{i,j} x_{i,j}$, increases while the constraints are still satisfied. However, the number of self-loops has decreased, contradicting the minimality of self-loops in $Q$. 
\end{proof}

\begin{claim}\label{even cycle}
For any $Q\in\mathcal{Q}$, there is no even cycle in $Q$. 
\end{claim}

\begin{proof}
Assume for contradiction that there exists $Q\in\mathcal{Q}$ containing an even cycle. Say $C = x_{1}x_{2} \dots x_{m-1}x_{m}x_1$ is an even cycle in $Q$. Let \[ s = p_{1,2}x_{1,2} + p_{2,3}x_{2,3} + \dots + p_{m-1,m}x_{m-1,m} + p_{1,m}x_{1,m} \] be the part of the objective function of \texttt{LP2-dual$'$}, $\sum p_{i,j} x_{i,j}$, involving this even cycle. Without loss of generality, assume $p_{1,2} + p_{3,4} + \dots + p_{m-3,m-2} + p_{m-1,m} \geq p_{2,3} + p_{4,5} + \dots + p_{m-2,m-1} + p_{1,m}$. Set $\varepsilon = \min(x_{2,3}, x_{4,5}, \dots, x_{m-2,m-1}, x_{1,m}) > 0$. Then
\begin{equation*}
    p_{1,2}(x_{1,2}+\varepsilon) + p_{2,3}(x_{2,3} - \varepsilon) + \dots + p_{m-1,m}(x_{m-1,m} + \varepsilon) + p_{1,m}(x_{1,m} - \varepsilon) \geq s.
\end{equation*}
Notice that the constraints are still satisfied. 
If this inequality is strict, then the optimality of $Q$ is contradicted. If there is equality, then, by our choice of $\varepsilon$, at least one of $x_{2,3}-\varepsilon,x_{4,5}-\varepsilon,\ldots,x_{1,m}-\varepsilon$ is 0. Thus we have obtained an optimal solution by removing an edge from $Q$, contradicting the minimality of non-loop edges in $Q$.
\end{proof}

Now our goal is to show that every connected component of $Q\in\mathcal{Q}$ with at least three vertices is an induced odd cycle. In the following claim, the degree of a vertex includes the possibility of self-loops.

\begin{claim}\label{pendant vertex}
For any $Q \in \mathcal{Q}$, if $C$ is a connected component of $Q$ with at least three vertices, then $C$ contains an odd cycle.
\end{claim}

\begin{proof}
Suppose that $C$ is a connected component of $Q\in\mathcal{Q}$ with at least three vertices.
We first prove that $C$ has no vertex of degree 1 in $Q$. Assume, to the contrary, that $C$ has a vertex $x_i$ of degree 1 in $Q$. Let $x_k$ be the neighbor of $x_i$. Then since $C$ is a connected component on at least three vertices, it follows that $x_k\in C$ and $x_k\neq x_i$. Thus $x_i$ does not have a self-loop, and so, by the construction of $Q$, we have $y_i = 0$. Thus the $i$th constraint $y_i + \sum_{i \neq j} x_{i,j} = 1$ in \texttt{LP2-dual$'$} becomes $x_{i,k} = 1$, and so the $k$th constraint in \texttt{LP2-dual$'$} implies that $x_k$ has degree 1 in $Q$. Hence, $C$ must be a connected component on only the two vertices $x_i,x_k$, a contradiction since $C$ has at least three vertices.
\par We now prove the claim. By Claim~\ref{even cycle}, it suffices to prove that $C$ contains a cycle.
Assume, to the contrary, that $C$ is acyclic. 
Then since $C$ is connected, acyclic, and $|C|\geq 2$, there are at least two vertices $x_{i_1}\neq x_{j_1}$ in $C$ of non-loop degree 1 in $C$. 
By the preceding paragraph, $x_{i_1}$ must have another neighbor $x_{i_2}$, and $x_{j_1}$ must have another neighbor $x_{j_2}$.
Since $C$ is a connected component, these two neighbors must lie in $C$, and so we must have
$x_{i_2}=x_{i_1}$ and $x_{j_2}=x_{j_1}$.
In other words, $x_{i_1},x_{j_1}$ have self-loops, contradicting Claim~\ref{self-loop}.
\end{proof}

\begin{claim}\label{claim:induced cycle}
For any $Q \in \mathcal{Q}$, every odd cycle in $Q$ is an induced cycle.
\end{claim}
\begin{proof}
Let $C$ be an odd cycle in $Q\in\mathcal{Q}$. Assume for contradiction that $C$ is not induced. Then there exists a chord in $E(C)$. This chord splits $C$ into two cycles $C'$ and $C''$. One of them is an even cycle, contradicting Claim~\ref{even cycle}. 
\end{proof}

In the remaining claims, we use the following notation. Define an edge weight function $w:V(Q)^2\to\mathbb{R}$ such that for each distinct $x_i,x_j\in V(Q)$, we have $w(x_i, x_j) = w(x_j, x_i) = x_{i,j}$ and $w(x_i, x_i) = y_i$. For convenience, we write $w(e)$ for $e\in E(Q)$.

\begin{claim}\label{disjoint}
For any $Q\in\mathcal{Q}$, the odd cycles in $Q$ are vertex disjoint.
\end{claim}

\begin{proof}
It is easy to see that if there are two distinct odd cycles in $Q$ sharing an edge or at least two vertices, then $Q$ contains an even cycle, contradicting Claim~\ref{even cycle}.

Let $C_1 = u_1u_{2}\cdots u_{2s+1}u_1$ and $C_2 = v_{1}v_{2}\cdots v_{2s'+1}v_{1}$ be distinct odd cycles in $Q\in\mathcal{Q}$. Assume that $C_1$ and $C_2$ share exactly one vertex, say $u_1 = v_1$. Let $\varepsilon = \min\{w(e): e\in E(C_1)\cup E(C_2)\}$. We show that by adjusting the values of the $w(e)$'s, we can either construct an optimal solution with fewer edges or increase the value of the objective function. In both cases, we reach a contradiction. 
We define a new edge weight function $w'$ as follows (where the indices of $u_{i+1},v_{i+1}$ below are considered modulo $2s+1, 2s'+1$ respectively):
\[ w'(e)= \begin{cases}
    w(e) + (-1)^{i}\varepsilon &\text{if }e=u_iu_{i+1}\in E(C_1),\\
    w(e) + (-1)^{i+1}\varepsilon &\text{if }e=v_iv_{i+1}\in E(C_2),\\
    w(e) &\text{otherwise}.
   \end{cases} \]
By choice of $\varepsilon$, we know that $w'(e) \geq 0$ for all $e\in E(Q)$. Furthermore, for each fixed $x_r\in V(Q)$, it is easy to see that $\sum_{r\neq i} w'(x_r,x_i) = \sum_{r\neq i} w(x_r,x_i)$. This shows that the constraints of $\texttt{LP2-dual}'$ are satisfied. Let $\Delta = \sum_{i< j}p_{i,j}(w'(x_i,x_j) - w(x_i,x_j))$. If $\Delta > 0$, then the optimality of $Q$ is contradicted. If $\Delta < 0$, then add 1 to the powers of $-1$ in the definition of $w'$ to contradict the optimality of $Q$. Assume $\Delta = 0$. By choice of $\varepsilon$, either there exists a $w'(e) = 0$, or there exists such a zero edge weight after adding 1 to the powers of $-1$ in the definition of $w'$. Since $\Delta=0$, we get an optimal solution with fewer non-loop edges, contradicting the minimality of non-loop edges in $Q$.
\end{proof}

\begin{claim}\label{odd cycle}
For any $Q\in\mathcal{Q}$, each connected component of $Q$ contains at most one odd cycle.
\end{claim}
\begin{proof}
Assume, for the sake of contradiction, that $C$ is a connected component of $Q\in\mathcal{Q}$ that contains distinct odd cycles $C_1$ and $C_2$. By Claim~\ref{disjoint}, these two cycles are vertex disjoint. Let $P$ be a shortest $(C_1, C_2)$-path. Let $C_1 = u_1u_{2}\cdots u_{2s+1}u_1, P = w_1w_{2}\cdots w_{\ell}$, and $C_2 = v_{1}v_{2}\cdots v_{2s'+1}v_{1}.$ Let $\delta = \min\{w(e): e\in E(C_1)\cup E(C_2)\}$ and $\delta' = \min\{w(e): e\in E(P)\}$. Define $\varepsilon = \delta$ if $2\delta \leq \delta'$, and $\varepsilon = \delta'/2$ otherwise. Without loss of generality, say $w_1 = u_1$ and $w_{\ell} = v_1$. We define a new edge weight function $w'$ as follows (where the indices of $u_{i+1},v_{i+1}$ below are considered modulo $2s+1, 2s'+1$ respectively):
\[ w'(e)= \begin{cases}
    w(e) + (-1)^{i}\varepsilon &\text{if }e=u_iu_{i+1}\in E(C_1),\\
    w(e) + (-1)^{\ell + i}\varepsilon &\text{if }e=v_iv_{i+1}\in E(C_2),\\
    w(e) + (-1)^{i+1}2\varepsilon &\text{if }e=w_iw_{i+1}\in E(P),\\
    w(e) &\text{otherwise}.
\end{cases}\]
We chose a shortest path $P$ to ensure that $P$ does not share edges with $C_1$ or $C_2$. This makes $w'$ well-defined. By choice of $\varepsilon$, we know that $w'(e) \geq 0$ for all $e\in E(Q)$. Furthermore, for each fixed $x_r\in V(Q)$, it is easy to see that $\sum_{r\neq i} w'(x_r,x_i) = \sum_{r\neq i} w(x_r,x_i)$. This shows that the constraints of \texttt{LP2-dual$'$} are satisfied. Let $\Delta = \sum_{i< j}p_{i,j}(w'(x_i,x_j) - w(x_i,x_j))$. If $\Delta > 0$, then the optimality of $Q$ is contradicted. If $\Delta < 0$, then add 1 to the powers of $-1$ in the definition of $w'$ to contradict the optimality of $Q$. Assume $\Delta = 0$. By choice of $\varepsilon$, either there exists a $w'(e) = 0$, or there exists such a zero edge weight after adding 1 to the powers of $-1$ in the definition of $w'$. Since $\Delta=0$, we get an optimal solution with fewer non-loop edges, contradicting the minimality of non-loop edges in $Q$.
\end{proof}

In the following claim, a pendant path $P$ of an odd cycle $C'$ is a path $v_1v_2\dots v_\ell$ such that $v_1$ is the only vertex on $P$ that lies on $C'$ and $v_\ell$ has degree 1 in non-loop edges.

\begin{claim}\label{no pendant path}
For any $Q\in\mathcal{Q}$, no connected component in $Q$ contains an odd cycle with a pendant path as a subgraph.
\end{claim}
\begin{proof}
Assume, for the sake of contradiction, that $C$ is a connected component of $Q\in\mathcal{Q}$ containing an odd cycle $C'$ with a pendant path $P$. Say $C' = u_1u_{2}\cdots u_{2s+1}u_1$ and $P = v_1 v_2\cdots v_{\ell}$ with $u_1 = v_1$. Let $\delta = \min\{w(e): e\in E(C')\}$ and $\delta' = \min\{\{w(e): e\in E(P)\}\cup \{w(v_\ell,v_\ell)\}\}$. We know that $w(v_\ell, v_\ell) > 0$ because $P$ is a pendant path. Define $\varepsilon = \delta$ if $2\delta \leq \delta'$, and $\varepsilon = \delta'/2$ otherwise. We define a new edge weight function $w'$ as follows (where the index of $u_{i+1}$ below is considered modulo $2s+1$):
\[ w'(e) = \begin{cases}
    w(e) + (-1)^{i}\varepsilon &\text{if } e=u_iu_{i+1}\in C',\\
    w(e) + (-1)^{i+1}2\varepsilon &\text{if } e=v_iv_{i+1}\in P,\\
    w(e) + (-1)^{\ell+1}2\varepsilon &\text{if } e=v_\ell v_\ell,\\
    w(e) &\text{otherwise}.
\end{cases} \]
By choice of $\varepsilon$, we know that $w'(e) \geq 0$ for all $e\in E(Q)$. Furthermore, for each fixed $x_r\in V(Q)\setminus v_{\ell}$, it is easy to see that $\sum_{r\neq i} w'(x_r,x_i) = \sum_{r\neq i} w(x_r,x_i)$. For $v_\ell$, we have $w'(v_{\ell-1}, v_{\ell}) + w'(v_{\ell}, v_{\ell}) = w(v_{\ell-1}, v_{\ell}) + w(v_{\ell}, v_{\ell})$. Thus the constraints of \texttt{LP2-dual$'$} are satisfied. Let $\Delta = \sum_{i<j}p_{i,j}(w'(x_i,x_j) - w(x_i,x_j))$. If $\Delta > 0$, then the optimality of $Q$ is contradicted. If $\Delta < 0$, then add 1 to the powers of $-1$ in the definition of $w'$ to contradict the optimality of $Q$. Assume $\Delta = 0$. By choice of $\varepsilon$, either there exists a $w'(e) = 0$, or there exists such a zero edge weight after adding 1 to the powers of $-1$ in the definition of $w'$. 
Since $\Delta=0$, we get an optimal solution with fewer edges, contradicting the minimality of edges in $Q$.
\end{proof}

\par We now have almost all the pieces needed to prove Claim~\ref{claim:Q'}. Fix $Q\in\mathcal{Q}$. Suppose that $C$ is a connected component of $Q$. If $C$ contains exactly one vertex, then $C$ is an isolated vertex with a self-loop. If $C$ contains exactly two vertices, then $C$ is an isolated edge which, a priori, may have self-loops. Suppose now that $C$ contains at least three vertices. By Claim~\ref{pendant vertex}, $C$ contains an odd cycle. By Claims \ref{claim:induced cycle} and \ref{odd cycle}, this cycle in $C$ must be induced and unique. By Claim~\ref{no pendant path}, every vertex in $C$ lies on this unique cycle. This implies that $C$ is an odd cycle which may have self-loops. In fact, the following claim proves that isolated edges and connected components with at least three vertices must have no self-loops.

\begin{claim}\label{self-loop 2}
For any $Q\in\mathcal{Q}$, if a vertex $v\in V(Q)$ has a self-loop, then $v$ is an isolated vertex.
\end{claim}

\begin{proof}
Let $Q\in\mathcal{Q}$, and let $v\in V(Q)$ be a vertex with a self-loop. Assume, to the contrary, that $v$ is not an isolated vertex. Then $v$ must lie in a connected component $C$ of $Q$ on exactly two vertices or at least three vertices. In the first case, we can use the constraints to deduce that both vertices in $C$ must have a self-loop. This contradicts Claim~\ref{self-loop}. Now assume that $C$ contains at least three vertices and $v\in C$. By the above discussion, $C$ is an odd cycle which may have self-loops.
Say $C=v_1v_2\cdots v_{2s+1}v_1$ for some positive integer $s$ with $v=v_1$. Let $\delta=\min\{w(e) : e\in E(C)\}$ and $\delta'=w(v,v)$. Define $\varepsilon = \delta$ if $2\delta \leq \delta'$, and $\varepsilon = \delta'/2$ otherwise. We define a new edge weight function $w'$ as follows (where the index of $v_{i+1}$ below is considered modulo $2s+1$):
\[ w'(e) = \begin{cases}
    w(e) + (-1)^i\varepsilon &\text{if } e=v_iv_{i+1}\in C,\\
    w(e) + 2\varepsilon &\text{if } e=vv,\\
    w(e) &\text{otherwise}.
\end{cases} \]
Now the proof is similar to that of Claim~\ref{no pendant path}. We can either increase the objective function or reduce the number of edges, a contradiction. 
\end{proof}

By Claim~\ref{self-loop 2}, only isolated vertices in $Q$ have self-loops. Constructing $Q'$ from $Q$, we get that $Q'$ is an odd cycle decomposition of $K_k$, thus proving Claim~\ref{claim:Q'}, as desired. This completes the proof of Proposition~\ref{asymptotic sharpness}, finishing the proof of Theorem~\ref{thm:main}.

\end{proof}

\noindent
{\bf Concluding Remarks:}
Observe that the $o(n_1\cdots n_k)$-term in the bound of Proposition~\ref{asymptotic sharpness} occurs due to divisibility issues. In particular, it is nonzero when an optimal solution $(a_1, \ldots, a_k)$ of problem \texttt{P} is not rational. Here are two examples where the bound is sharp for infinitely many values of $n_1,\ldots, n_k$. 
\begin{enumerate}
    \item Suppose that, in the hypotheses of Theorem~\ref{thm:main}, $k\geq 3$ is an odd integer, and for all $i\neq j\in[k]$, $1-d_{i,j} = a^2/b^2 \in\mathbb{Q}$. That is, all the $1-d_{i,j}$'s are the same rational perfect square.
    By Theorem~\ref{thm:main}, the number of independent transversals is at most $\prod_{i=1}^k \sqrt{1-d_{i, i+1}}\prod_{i=1}^k n_i = (a/b)^k\prod_{i=1}^k n_i$, where the index $i+1$ is modulo $k$. We now construct a graph $G$ that attains this bound. Note that there are infinitely many integers $n$ such that $(a/b)n$ is an integer. Let $n$ be such an integer, and for all $i\in[k]$, let $n_i = n$ and $S_i\subseteq V_i$ be a set of size $(a/b)n$. For each $i\neq j\in[k]$, put a complete bipartite graph between $V_i\setminus S_i, V_j$ and between $S_i,V_j\setminus S_j$ with no other edges occurring between $V_i,V_j$. Then $G$ satisfies $1-d_{i,j}=a^2/b^2$ for all $i\neq j\in[k]$, and the number of independent transversals in $G$ is $\prod_{i=1}^k|S_i| = (a/b)^kn^k$, showing that $G$ attains the desired bound.
    \item Suppose the hypotheses of Theorem~\ref{thm:main}. Let $I_1,\ldots, I_\ell$ be a a partition of $[k]$ such that for all $r\in[\ell]$, $|I_r|$ is odd. For all $r\in[\ell]$ and $i\neq j\in I_r$, let $1-d_{i,j} = a^2_r/b^2_r\in\mathbb{Q}$ (if $|I_r|=1$, let $a_r/b_r=1$); otherwise, if $i\in I_r, j\not\in I_r$, let $1-d_{i,j} = 1$. We are essentially taking disjoint copies of example~1 above. By Theorem~\ref{thm:main}, the number of independent transversals is at most $\prod_{r=1}^{\ell}(a_r/b_r)^{|I_r|}\prod_{i=1}^k n_i$. We construct a graph $G$ that attains this bound for infinitely many $n_1,\ldots,n_k$ by applying example~1 to each $I_r$. 
\end{enumerate}

One could also study stability type of results. One could show that if $d_{i,j} = d$ for all $i,j$, then, by an argument similar to our proof, if the number of independent transversals in $G$ is close to the optimal value, then $G$ should be close to the construction in Section~\ref{subsec:as} where $a_i = \sqrt{1-d}$ for each $i$. In the case where the $d_{i,j}$ are not all equal to each other, there could be very different optimal constructions. However, one should be able to show that each of those optimal constructions is close to a construction corresponding to some solution of \texttt{P}.

\noindent
{\bf Acknowledgment:}
This project was done as part of the 2021 New York Discrete Math REU, funded by NSF grant DMS 2051026. The authors would like to thank Drs. Adam Sheffer and Pablo Sober\'{o}n for organizing the REU. The authors would also like to thank the anonymous referees for their useful comments and, in particular, Reviewer B for a comment that helped simplify the proof of Lemma~\ref{lem:cycle bound}. 

\pagebreak

\end{document}